\documentclass[a4paper,11pt]{article}
\usepackage{amssymb}

\usepackage{graphicx}

\newcommand{\C}{\mathbb{C}}
\newcommand{\R}{\mathbb{R}}

\newcommand{\N}{\mathbb{N}}

\newcommand{\cuad}{{\sqcap\kern-.68em\sqcup}}

\newcommand{\norm}[1]{\|#1\|}

\newcommand{\cJ}{{\mathcal J}}

\newtheorem{theorem}{Theorem}[section]
\newtheorem{proposition}{Proposition}[section]

\newtheorem{lemma}{Lemma}[section]

\newtheorem{remark}{Remark}[section]
\newcommand{\bremark}{\begin{remark} \em}
\newcommand{\eremark}{\end{remark} }

\def\JJ{\mathbb{J}}

\usepackage{cite}
\begin{document}

\begin{center}{\bf  \Large   On isolated  singular solutions of semilinear\smallskip

 Helmholtz equation

 }\bigskip

 \bigskip

 {\small

{\sc  Huyuan Chen\footnote{chenhuyuan@yeah.net}}

\medskip
Department of Mathematics, Jiangxi Normal University,\\
Nanchang, Jiangxi 330022, PR China \\[10pt]
 
 {\sc  Feng Zhou\footnote{fzhou@math.ecnu.edu.cn} }

\medskip
Center for PDEs, School of Mathematical Sciences, East China Normal University,\\
Shanghai Key Laboratory of PMMP, Shanghai 200062, PR China \\[20pt]
}

\begin{abstract}
 
Our purpose of this paper is to study   isolated singular solutions  of semilinear Helmholtz equation
$$ -\Delta 
u-u=Q|u|^{p-1}u \quad{\rm in}\ \ \R^N\setminus\{0\},\ \qquad\lim_{|x|\to0}u(x)=+\infty, $$
where $N\geq 2$, $p>1$ and  the potential $Q:\R^N\to (0,+\infty)$ is a H\"older continuous function satisfying extra decaying conditions at infinity. We give the classification of the   isolated singularity in the Serrin's subcritical case and then  isolated singular solutions is derived with the form $u_k=k\Phi+v_k$ via the Schauder fixed point theorem for  the integral equation
$$v_k=\Phi\ast\big(Q|kw_\sigma+v_k|^{p-1}(kw_\sigma+v_k)\big)\quad{\rm in}\ \, \R^N,$$ 
 where $\Phi$ is the real valued fundamental solution $-\Delta-1$ and $w_\sigma$ is a also a real valued solution $(-\Delta-1)w_\sigma=\delta_0$ with the asymptotic behavior at infinity controlled by $|x|^{-\sigma}$ for some $\sigma\leq \frac{N-1}{2}$.

 \end{abstract}

  \end{center}
  
  \vspace{1mm}


\vspace{2mm}

\setcounter{equation}{0}
\section{Introduction}
Our purpose of this paper is to study   isolated singular solutions  of
semlinear Helmholtz equation 
 \begin{equation}\label{eq 1.1}
\left\{
 \arraycolsep=1pt
\begin{array}{lll}
\displaystyle \ \ \ -\Delta u-u=Q|u|^{p-1}u \quad   \rm{in}\ \, \R^N\setminus\{0\},\\[2mm]
 \phantom{     }
\displaystyle
\lim_{|x|\to0^+}u(x)=+\infty,
\end{array}\right.
 \end{equation}
where $N\geq 2$, $p>1$ and the potential $Q:\R^N \to \R$ is a H\"older continuous function.

The isolated singularities of semilinear elliptic equations have been studied extensively in the last decades since the first study  due to Brezis in an unpublished
 note  (see the introduction in \cite{BB}), motivated by  Thomas-Fermi problem \cite{Br},
 where proves that the  elliptic problem with absorption nonlinearity
$$
 -\Delta  u+|u|^{q-1}u=k\delta_0 \quad   \rm{in}\ \, \Omega,\qquad
u=0  \quad   \rm{on}\ \, \partial\Omega 
$$
admits a unique distributional  solution $u_k$ for $1< q< N/(N-2)$ and $k>0$, while no solution exists when $q\geq  N/(N-2)$, where $\delta_0$ is the Dirac mass concentrated at the origin. {\it Here $N/(N-2)$ is the Serrin's critical exponent. } More analysis on differential equations with measures could  see \cite{KR,PW,V1} and \cite{V} for a survey. Isolated singularity of problem
\begin{equation}\label{local}
 -\Delta  u+|u|^{q-1}u=0 \ \   \rm{in}\ \, \Omega\setminus\{0\}, 
 \qquad
u=0  \ \   \rm{on}\ \, \partial\Omega
\end{equation}
has been classified by Brezis and V\'{e}ron  in \cite{BV},  that  problem (\ref{local})
 admits only the zero solution when $q\geq  N/(N-2)$. When $1< q< N/(N-2)$, V\'{e}ron in \cite{V0} described   all the possible singular behaviors of positive  solutions of (\ref{local}). In particular he proved that the singularity must be isotropic when
$(N+1)/(N-1)\leq q< N/(N-2)$, and under the positivity assumption  two types of singular behaviors occur:\smallskip

\noindent (i) either $u(x)\sim c_Nk|x|^{2-N}$ as $|x|\to 0$ and $k$ can take any positive value; $u$ is said to have a {\it weak singularity} at $0$, and actually $u=u_k$,

\noindent (ii) or $u(x)\sim c_{N,q}|x|^{-\frac2{q-1}}$ as $|x|\to 0$;  $u$ is said to have a {\it strong singularity} at $0$, and
$\displaystyle u= \lim_{k\to\infty}u_k$.\\

In contract with absorption nonlinearity,  the elliptic problem with source nonlinearity 
\begin{equation}\label{eq 1.4}
  -\Delta   u=u^p\ \   {\rm in}\ \, \Omega\setminus\{0\}\qquad 
   u= 0 \ \   {\rm on}\ \, \partial\Omega 
\end{equation}
has a different structure of   the isolated singular solutions.  It was  classified  by Lions in \cite{L}   by building the connections with the weak solutions of
\begin{equation}\label{eq 1.5}
 \displaystyle    -\Delta    u=u^p+k\delta_0\ \
 {\rm in}\ \, \Omega,\qquad u= 0 \ \  {\rm on}\ \, \partial\Omega.
\end{equation}
For $N\geq3$ and $p\in (1,\frac{N}{N-2})$  any positive  solution of (\ref{eq 1.4})
is a weak solution of (\ref{eq 1.5}) for some $k\ge0$; when $p\ge \frac{N}{N-2}$, the parameter $k=0$. In the subcritical $p\in(1,\frac{N}{N-2})$, the solution of (\ref{eq 1.4}) has either the singularity of $|x|^{2-N}$ or   removable singularity. For the existence,  Lions in \cite{L} showed that  there exists $k^*>0$ such that,  for $k\in(0,k^*)$ problem (\ref{eq 1.5}) has at least two positive  solutions including the minimal solution and a Mountain Pass type solution; for $k=k^*$  problem (\ref{eq 1.5}) has a unique solution; for $k>k^*$ there is no solution of  (\ref{eq 1.5}). Setting $\Omega=B_r(0)$ and $k^*(r)$ the critical number ensuring the existence of positive solution of (\ref{eq 1.5}), it is well-known that 
\begin{equation}\label{re 1.1--r}
\lim_{r\to+\infty}k^*(r)=0.
\end{equation}
 In the super critical case,   the 
 singularity at the origin of positive solutions to
 \begin{equation}\label{eq 1.5-super}
   -\Delta    u=u^p \ \ {\rm in}\ \, B_1(0)\setminus\{0\}
\end{equation}
 is classified for  $p = \frac{N}{N-2}$ by Aviles in \cite{Av},
 for $\frac{N}{N-2} < p < \frac{N+2}{N-2}$ by  Gidas and  Spruck in \cite{GS}
 and for $p=\frac{N+2}{N-2}$ by Caffarelli, Gidas and  Spruck in \cite{CGS}. \smallskip

    The Helmholtz equations  arise from  the nonlinear wave equation
$$
\frac{\partial^2 }{\partial t^2}\psi(t,x) - \Delta \psi(t,x)  =h(x,|\psi|^2)\psi, \qquad (t,x)\in\R\times\R^N,
$$
which is used to describe wave propagation in an ambient medium with nonlinear response,
where $h:\R^N\times \R\to\R$ is a real-valued function. The time-periodic ansatz 
$\psi(t,x)=e^{-i \kappa t}u(x)$ with $\kappa >0$
leads to the nonlinear Helmholtz equation.
 Motivated  their great applications,  the nonlinear Helmholtz equations  have been attracted attentions by researchers,  see \cite{A,CEW,EP0,E,EW0,EW1,EY,G,J,M1,MMP}
 and the references therein. In particular, \cite{EW1} applies the dual variational methods  to obtain the weak real-valued solution of
 the integral equation of 
  $-\Delta u-k^2 u= Q|u|^{p-1}u$ in $\R^N$ to avoid
 the strong indefiniteness of the normal energy functional works on $H^1(\R^N)$. For 
  $p\in(\frac{N+3}{N-1}, \frac{N+2}{N-2})$, a weak solution of 
 $$u=\Phi\ast  (Q |u|^{p-1}u), $$ 
 letting $v=Q^{\frac1{p+1}} |u|^{p-1}u$,  by obtaining the critical point of the energy functional  
$$
\cJ_{t}(v) =\frac{1}{q'}\int_{\R^N}|v|^{q'}\, dx + \frac12 \int_{\R^N}(Q^{\frac1{p+1}}v) \Phi\ast (Q^{\frac1{p+1}} v)\, dx,\quad\forall\, v\in L^{q'}(\R^N), $$
 where $q'=\frac{p+1}{p}$ is the conjugate index of $p+1$,  $\Phi$ is a radially symmetric fundamental solution of Helmholtz operator in $\R^N$ verifying 
 \begin{equation}\label{funda 1.1}
  \Phi(x)\sim \left \{
  \begin{array}{lll}
 c_N |x|^{2-N} &\ \  {\rm if} \ \ N\geq 3\\[2mm]
 -c_N\ln|x|    &\ \  {\rm if} \ \ N=2
  \end{array}
\right. \ \ {\rm as}\ \ |x|\to0^+, \qquad \limsup_{|x|\to+\infty}|\Phi(x)||x|^{\frac{N-1}{2}}<+\infty
 \end{equation}
with the normalized constant $c_N>0$.  
 More properties of the fundamental solution is mentioned in Section 2. 
  In a recent paper \cite{CEW}, the authors studied the complex valued weak solution of  $-\Delta u-k^2 u= Qf(x,u)$ in $\R^N$, where $f$ could involve the perturbation of non-homogeneous source with the model
 $$f(x,u)=|u|^{p-1}u+g(x)$$
  with  $g$ be bounded and decaying at infinity.

Thanks to (\ref{re 1.1--r}), there is no hope to find positive solutions of Helmholtz equation (\ref{eq 1.1}). Our aim in this article is to obtain sign-changing solutions with isotropic singularity at the origin in the Serrin's subcritical case, i.e. $p\in(1,p^*_{_N})$, 
where  the Serrin's exponent  $p^*_{_N}$ is defined by
 \begin{equation}\label{serrin 1} 
 p^*_{_N}=
\arraycolsep=1pt\left\{
\begin{array}{lll}
  \frac{N}{N-2} \quad
 &{\rm if}\ \,  N\geq 3,\\[2mm]
 \phantom{  }
 \displaystyle   +\infty \quad
 &{\rm if}\ \, N=2.
\end{array}\right.
 \end{equation}
Our strategy of deriving isolated singular solution of (\ref{eq 1.1}) is to build the connection with  weak solutions of 
\begin{equation}\label{eq 1.1w}
  -\Delta   u-u =Q|u|^{p-1}u+k\delta_0\quad \rm{in}\ \, \R^N,
 \end{equation}
where $k>0$ and $\delta_0$ is the Dirac mass at the origin. Here we say that  $u\in L^1_{loc}(\R^N)$ is  
a weak solution of (\ref{eq 1.1w}) if 
$$
 \int_{\R^N} \left[u(-\Delta) \xi-u\xi -Q |u|^{p-1}u\xi\right]\,dx=  k \xi(0),\quad  \forall\, \xi\in C^\infty_c (\R^N).
$$

Our first result is on isotropic singular solutions of  (\ref{eq 1.1}) classified by 
Dirac source in the distributional sense.

 \begin{theorem}\label{teo 1}
Assume that $N\ge2$,  $p\in(1,\, p^*_{_N})$,  $Q:\R^N \to \R$ is a  H\"older continuous function in $\R^N$ such that $Q(0)>0$.

If $u$ is a  classical solution of (\ref{eq 1.1}), then  $u\in L^1_{loc}(\R^N)$, $Q|u|^{p-1}u\in   L^1_{loc}(\R^N)$ and  $u$ is a very weak solution of (\ref{eq 1.1w}) for some $k>0$ with the asymptotic behaviors that 
     $$\lim_{|x|\to0^+} \frac{u(x)}{\Phi(x)} = k,$$
where $c_N>0$ is a normalized constant.

\end{theorem}
 
 \begin{remark}
We remark that the regularity of $Q$ could be weakened to be $L^\infty$, 
 Moreover, the asymptotic behavior at the origin could be replaced by 
 $$Q(x)\geq Q_0|x|^{\varrho_0}\quad {\rm in}\ \, B_{r}(0)\setminus\{0\} $$
for some $Q_0>0$, $r>0$ and $\varrho_0>-2$.  In this case, the Serrin's exponent $p^*_N$ should be
  $$ 
  \frac{N+\varrho_0}{N-2} \quad {\rm if}\ \,  N\geq 3,\qquad
 +\infty \quad
  {\rm if}\ \, N=2.$$
For simplicity, we in this paper deal with the case that  $Q$ is continuous and $Q(0)>0$. 
\end{remark}
   
  From the Fourier transform, the Helmholtz harmonic functions, {\it i.e. $u$ is Helmholtz harmonic in $\R^N$ if}
\begin{equation}\label{1.2}
-\Delta u-u=0\quad{\rm in}\ \ \R^N,
\end{equation}
 have a general expression 
 $u(x)=\int_{\mathbb{S}^{N-1}}e^{i\xi\cdot x} g(\xi) d\omega(\xi) $
 by the property that 
  ${\rm supp}(\mathcal{F}(u))\subset \mathbb{S}^{N-1}$, where $\mathbb{S}^{N-1}$ is the unit sphere in $\R^N$, see \cite{A}.
  To this end, an appreciated  classification of  the Helmholtz harmonic function is needed  from the behaviors at infinity and we denote ${\mathcal S}_\sigma$ with $\sigma\in\R$,  the set of all nontrivial solutions of homogeneous Helmholtz equation (\ref{1.2}) with the asymptotic behavior controlled by $(1+|x|)^{-\sigma}$, i.e.
 \begin{equation}\label{sol set}
 \mathcal{S}_\sigma=\Big\{u\  {\rm is\ Helmholtz\ harmonic}:\, \sup_{x\in\R^N}\left(|u(x)|(1+|x|)^{ \sigma}\right)<+\infty\Big\}.
 \end{equation}
 \begin{proposition}\label{pr hom sol}

 Let $N\geq 2$ and ${\mathcal S}_\sigma$ be the solutions' set defined by (\ref{sol set}),
 then
 $$\mathcal{S}_\sigma=\emptyset\quad {\rm for}\ \,\sigma>\frac{N-1}{2}$$
and
 the mapping: $\sigma\in(-\infty,\frac{N-1}{2}]\to {\mathcal S}_\sigma$ is strictly decreasing, i.e. for $-\infty<\sigma_2<\sigma_1< \frac{N-1}{2}$, we have
that
$$ \emptyset \not= \mathcal{S}_{\frac{N-1}{2}}\varsubsetneqq {\mathcal S}_{\sigma_1}\varsubsetneqq{\mathcal S}_{\sigma_2}.$$
 \end{proposition}

 \begin{remark}
Note that the typical non-vanishing Helmholtz functions are
$\psi(x)=\sin x_i$, $\psi(x)=\cos x_i$ for $x=(x_1,\cdots, x_N)$, $i=1,\cdots,N$.
\end{remark}

 For the existence of weak solution of (\ref{eq 1.1w}), we need more restrictions on $Q$ as following:
 
 \begin{itemize}
\item[$( \mathbb{Q}_\alpha)$]  $Q:\R^N \to\R$ is a H\"older continuous function in $ L^\infty_\alpha(\R^N)$ with $\alpha\in \R$ such that 
$Q(0)>0$,
where   $L^\infty_\alpha(\R^N)$ is the space of functions $w$
such  that
 $$\norm{w}_{L^\infty_\alpha}:={\rm esssup}_{x\in\R^N}|w(x)|(1+|x|)^{\alpha}<+\infty.$$
\end{itemize}
For a given number $\alpha\in\R$, we denote  
 \begin{equation}\label{serrin 2} 
 p^\#_\alpha=1+\frac{2}{N-1}\Big(\frac{N+1}{2}-\alpha\Big). 
 \end{equation}
Note that the Serrin exponent $p^*_{_N}$ arises from the isolated singularity related to the fundamental solution at the origin and while  the critical exponent $p^\#_\alpha$ is the one to control the decay at infinity. For the existence of singular solution,  it requires  $p^\#_\alpha<p^*_{_N}$,  which  holds if
  \begin{equation}\label{req 3}
  \alpha\in\R\ \ {\rm for}\  N=2,\quad \alpha>0\ \ {\rm for}\  N=3 \quad{\rm and}\quad \alpha>\frac{N(N-3)}{2(N-2) }\ \ {\rm if}\ \ N\geq 4. 
  \end{equation}

 \begin{theorem}\label{teo 2}
Assume that $N\geq  2$,  the potential $Q$  verifies $( \mathbb{Q}_\alpha)$
 with $\alpha$ satisfying  (\ref{req 3}). 
 Let $p^\#_\alpha\, p^*_{_N}$ be defined in (\ref{serrin 1}) and (\ref{serrin 2}) respectively.
 
 Then for 
 $$p\in \big(1,\,+\infty\big)\cap \big[p^\#_\alpha,\, p^*_{_N} \big)$$
and   
$\psi_\sigma\in \mathcal{S}_\sigma$ 
with 
$$\sigma\in\left[ \frac{1}{p-1}\Big(\frac{N+1}{2} -\alpha\Big) ,\,  \frac{N-1}{2}\right],$$
 there exists $k^*>0$
  such that  for $ k\in(0, k^*)$, problem (\ref{eq 1.1w}) admits a  solution $u_{k}$
  in a integral form
 \begin{equation}\label{eq integral}
 u_k=\Phi\ast (Q|u_k|^{p-2}u_k)+k(\Phi+\psi_\sigma)\quad{\rm in}\ \, \R^N. \end{equation}

Moreover, $(i)$  $u_{k}$ is a classical solution of (\ref{eq 1.1}) and satisfy
 \begin{equation}\label{1.2-e}
 \lim_{x\to0} \frac{u_k(x)}{\Phi(x)}= k 
 \quad{\rm and}\quad   
 \limsup_{x\to0}  |u_k(x)||x|^{\sigma}<+\infty;
 \end{equation}
$(ii)$ if $\sigma\in\left( \frac{1}{p-1}\Big(\frac{N+1}{2} -\alpha\Big) ,\,  \frac{N-1}{2}\right)$,
  \begin{equation}\label{1.2-f}
  \limsup_{|x|\to+\infty} |u_k(x)-k(\Phi+\psi_\sigma)(x)||x|^{\sigma_p}<+\infty,
   \end{equation}
 where  $ \sigma_p:= \min\left\{\frac{N-1}{2},\alpha+\sigma p-\frac{N+1}2\right\}>\sigma$.

\end{theorem}

  \begin{remark}
 $(i)$   $u_k$ is disturbed from the singular function
 $kw_\sigma:=k(\Phi+\psi_\sigma)$, where $\psi_\sigma\in \mathcal{S}_\sigma$, $\Phi$ is the radially symmetric fundamental solution of Helmholtz operator and  $w_\sigma$ is also a weak solution of 
 $$-\Delta u-u=\delta_0\quad{\rm in}\ \, \R^N.$$

$(ii)$  If  $\alpha\geq \frac{N+1}{2}$, then for any $p\in(1,p^*_{_N})$
 the solution could be disturbed by functions in $\mathcal{S}_0$.  This means that problem (\ref{eq 1.1}) admits some singular solutions, which don't have to vanish uniformly at infinity. 
 
 Moreover, when $p>p^\#_\alpha$,
 our choice of $\sigma$ verifies that  $\sigma< \sigma_p$, so for $\psi_\sigma\in \mathcal{S}_{\sigma_p}\setminus  \mathcal{S}_{\sigma}$, the singular solutions could be differed by the asymptotic behaviors of $\psi_\sigma$. This means, if $\psi_{\sigma,1},\psi_{\sigma,2}\in\mathcal{S}_{\sigma_p}\setminus  \mathcal{S}_{\sigma}$ such that  $\psi_{\sigma,1}-\psi_{\sigma,2}\in\mathcal{S}_{\sigma_p}\setminus  \mathcal{S}_{\sigma}$, then corresponding solutions are different. 
 
 While for $p=p^\#_\alpha$ and $\sigma=\sigma_p$, our theorem shows the existence but we can not 
 distinguish them by the asymptotic behaviors directly.

  $(iii)$ When $N=2$, we can deal with the potential $Q\equiv 1$. 
\end{remark}

Our approach for singular solutions of (\ref{eq 1.1w}) is to use the Schauder fixed point theorem of (\ref{eq integral})  for the integral inequation
$$v_k=\Phi\ast\big(Q|kw_\sigma+v_k|^{p-2}(kw_\sigma+v_k)\big)\quad{\rm in}\ \, \R^N,$$ 
 where $w_\sigma= \Phi+\psi_\sigma $. This method also could be applied to obtain the complex valued solutions of (\ref{eq 1.1}), which is discussed in Section 5.\smallskip  
 
The rest of this paper is organized as follows. In Section 2, we  introduce the fundamental solution of Helmholtz operator and classify the the set of Helmholtz harmonic functions.     Section 3 is devoted to the classification of isolated singularity of classical solution of (\ref{eq 1.1}).   In Section 4, we  prove the existence of distributional solution of 
(\ref{eq 1.1w}) by Schauder fixed point theorem.
 
\setcounter{equation}{0}
\section{Preliminary}
\subsection{Fundamental solutions}

In polar coordinates, we have that
$$\Delta=\frac{\partial^2}{\partial r^2}+\frac{N-1}{r}\frac{\partial}{\partial r}+\frac{1}{r^2}\Delta_{\mathbb{S}^{N-1}}, $$
where
$\Delta_{\mathbb{S}^{N-1}}$ is the Beltrami Laplacian and $\mathbb{S}^{N-1}=\{x\in\R^N:\, |x|=1\}$. It is known that the eigenvalue $\mu_j=j(j+N-2)$ and corresponding eigenfunctions $v_{j}$ for $\Delta_{\mathbb{S}^{N-1}}$, i.e.
$$
 -\Delta_{\mathbb{S}^{N-1}} v_{j} =\mu_j  v_{j} \quad  {\rm in } \ \, \mathbb{S}^{N-1}.
 $$
 In fact, the multiplicity of eigenfunctions relative to $\mu_j$ equals the dimension of the space of homogenous, harmonic polynomials of degree $j$.

Let
\begin{equation}\label{z 1.1}
\mathcal{J}_{\Lambda_j}(t)=\sum^{+\infty}_{i=0}\frac{(-1)^i}{i!\Gamma(i+\Lambda_j+1)}\left(\frac t2\right)^{2i+\Lambda_j}
\end{equation}
which is the Bessel function of first kind with order
 \begin{equation}\label{1.1-1}
\Lambda_j:=\sqrt{\frac{(N-2)^2}{4}+\mu_j}\ ,
\end{equation}
that is,
$$
\frac{d^2}{d r^2}\mathcal{J}_{\Lambda_j}+\frac{1}{r}\frac{d}{d r}\mathcal{J}_{\Lambda_j}+(1-\frac{\Lambda_j^2}{r^2})\mathcal{J}_{\Lambda_j}=0.
$$

For $j\in\N$, denote
 \begin{equation}\label{1.1-2}
\psi_{j}(x)=|x|^{-\frac{N-2}{2}}\mathcal{J}_{\Lambda_j}(|x|)v_j(\frac{x}{|x|}), \quad\forall\, x\in\R^N,
 \end{equation}
 then  $\psi_{j}$ is a classical solution of Homogeneous Helmholtz equation (\ref{1.2}). Particularly, we set 
 \begin{equation}\label{1.1-2-0}
\Psi=\psi_{0}\quad {\rm in}\ \, \R^N.
 \end{equation}

It is known that
\begin{equation}\label{1.2-1}
 \lim_{r\to0^+}\psi_{j}(r)r^{\frac{N-2}{2}-\Lambda_j}=\frac{2^{-\Lambda_j}}{\Gamma(\Lambda_j+1)}
 \end{equation}
 and $\psi_{j}$ is oscillated at $+\infty$ with the order $r^{-\frac{N-1}{2}}$.

Let $\Lambda_N=\frac{N-2}{2}$ and denote 
$$\Phi(x)=c_0 |x|^{-\frac{N-2}{2}}  \lim_{\Lambda\to \Lambda_N}\frac{ \cos (\Lambda\pi) \JJ_{\Lambda}(|x|)-\JJ_{-\Lambda}(|x|)  }{\sin{\Lambda \pi}},\quad \forall \,x\in\R^N\setminus\{0\},$$
where $c_0>0$ depends on $N$,  $\mathcal{J}_{\frac{N-2}{2}}$ is defined by (\ref{z 1.1}) replaced $\Lambda_k$ by $\frac{N-2}{2}$ and
$$\mathcal{J}_{-\Lambda}(t)=\sum^{+\infty}_{m=0}\frac{(-1)^m}{m!\Gamma(m-\Lambda+1)}\left(\frac t2\right)^{2m-\Lambda},\quad \forall\, t>0.$$
Note that  the related asymptotic behaviors for
Bessel functions show that as $ |x|\to0^+$
\[
\Phi(x)=\left \{
  \begin{array}{lll}
 c_N |x|^{2-N}+O(|x|^{3-N}) &\quad   {\rm if} \ \, N\geq 3,\\[2mm]
 -c_N\ln|x|+O(|x|^{3-N})   &\quad  {\rm if} \ \, N=2,
  \end{array}
\right. 
\]
and 
$-\Delta \Phi-\Phi=\delta_0$ in the distributional sense. Moreover, we have $-\Delta \Phi-\Phi=0$ in $\R^N\setminus\{0\}$ pointwisely.
  As a conclusion,   $\Phi$ is a real fundamental solution for the
Helmholtz equation, so if $v$ is, say, a Schwartz function
on~$\R^N$ the convolution $\Phi*v$ satisfies
\begin{equation}\label{fundamental}
-\Delta(\Phi*v)-\Phi*v=v\quad{\rm in}\ \R^N.
\end{equation}
The function $\Phi$ is radially symmetric,  has the singularity $c_N|x|^{2-N}$ near the origin and oscillates at infinity controlled by
$|x|^{\frac{1-N}{2}}$.

 On the other hand, a complex-valued fundamental solution $\Phi_c$ could be involved as
 $$
\Phi_c(x):= (2\pi)^{-\frac{N}{2}} \mathcal{F}^{-1}((|\xi|^2-1-{\rm i}0)^{-1})(x)
=\frac{{\rm i}}{4} (2\pi |x|)^{\frac{2-N}{2}}H^{(1)}_{\frac{N-2}{2}}(|x|),\ \, \forall\,x \in \R^N \setminus \{0\},
$$
where i is the imagine unit,  $H^{(1)}_{\frac{N-2}{2}}$ is the Hankel function of the first kind of order $\frac{N-2}{2}$ and
$\mathcal{F}^{-1}((|\xi|^2-1-i0)^{-1})$ is approached by $[-\Delta-(1+{\rm i}\epsilon)]^{-1}$ as $\epsilon\to0^+$.
 In fact, the operator $-\Delta-(1+{\rm i}\epsilon)$: $H^2(\R^N)\subset L^2(\R^N)$
$\to$ $L^2(\R^N)$ is an isomorphism. Moreover, for any $f$ from the Schwartz space $\mathcal{S}$ its inverse is given by
$$
 [-\Delta-(1+{\rm i}\epsilon)]^{-1}f(x) = (2\pi)^{-\frac{N}{2}}
\int_{\R^N}e^{{\rm i} x\cdot \xi}\frac{\hat{f}(\xi)}{|\xi|^2-(1+{\rm i}\epsilon)}\, d\xi.
$$

For $H^{(1)}_{\frac{N-2}{2}}$, we have the asymptotic expansions
$$
  \label{eq:3}
H^{(1)}_{\frac{N-2}{2}}(s) =\left \{
  \begin{array}{lll}
 \sqrt{\frac{2}{\pi s}}\,
  e^{{\rm i} (s-\frac{N-1}{4}\pi)}[1+O(s^{-1})] & \quad  {\rm as} \ \ s \to
    \infty,\\[2.5mm]
 -\frac{{\rm i}\, \Gamma(\frac{N-2}{2})}{\pi}
\Bigl(\frac{2}{s}\Bigr)^{\frac{N-2}{2}}[1+O(s)]  & \quad  {\rm as} \ \ s \to 0^+,
  \end{array}
\right.
$$
(see e.g. \cite[Formulas (5.16.3)]{lebedev}), so there exists a constant $c_0>0$ such that
\begin{equation}
  \label{eq:4}
  |\Phi_c(x)| \le c_0 \max \left\{|x|^{2-N},1+(-\ln|x|)_+, (1+|x|)^{\frac{1-N}{2}}\right\} \quad
   {\rm for }\ x \in \R^N \setminus \{0\},
\end{equation}
where   $ t_+=\max\{0,  t\}$.
Moreover, $\Phi_c$ satisfies the equation $-\Delta \Phi_c - \Phi_c = \delta_0$
together with Sommerfeld's outgoing
radiation condition
\begin{equation}
  \label{eq:24}
|\nabla \Phi_c(x)- i \Phi_c(x)\hat x| = o(|x|^{\frac{1-N}{2}}) \quad
 {\rm as }\  |x| \to \infty,
\end{equation}
where $\hat x=\frac{x}{|x|}$. 
Furthermore, one has that 
$$\Phi_c=\Phi+i\Psi,$$
 where $\Psi$ is a real Helmholtz harmonic function such that,  by (\ref{eq:24}), 
 $$
|\nabla \Phi(x)+ \Psi(x)\hat x| = o(|x|^{\frac{1-N}{2}}), \quad |\nabla \Psi(x)-\Phi(x)\hat x| = o(|x|^{\frac{1-N}{2}}) \quad
 {\rm as}\  \, |x| \to \infty.
$$


\subsection{  Sets of  Helmholtz harmonic functions}
In this subsection, we will show the existence of the various  homogeneous Helmholtz functions. \medskip

\begin{lemma}\label{lm 2.1-1}
For $\sigma\in\R$, denote 
\begin{equation}\label{sol set1}
 \mathcal{S}^o_\sigma=\Big\{u\  {\rm is\ Helmholtz\ harmonic}:\, \limsup_{|x|\to+\infty}\left(|u(x)|(1+|x|)^{ \sigma}\right)\in(0,+\infty)\Big\}.
 \end{equation}
Then  for $-\infty<\sigma<  \frac{N-1}{2}$, we have
that
$$ \mathcal{S}^o_\sigma \not= \emptyset.$$  
 
\end{lemma}
\noindent{\bf Proof.}  Particular, for any positive integer $M\leq N$, we can set
$\Phi_M$ defined in (\ref{1.1-2-0}) replacing $N$ by $M$,  is also a Helmholtz harmonic in $\R^{N}$ with the decay at infinity controlled by   $|x|^{-\frac{M-1}{2}}$.

Generally, we next  show when $\sigma\in(-\frac{N-1}{2},+\infty)$, there is a Helmholtz harmonic function
$u$ such that
$$\limsup_{|x|\to+\infty} |u(x)|(1+|x|)^{-\sigma}\in(0,+\infty). $$

 Since
$$\Psi(x)=|x|^{-\frac{N-1}{2}}\sin(|x|-\varrho_0)+O(|x|^{-\frac{N+1}{2}})\quad{\rm near\ the\  infinity}, $$
where $\varrho_0=\frac{\Lambda_0\pi}{2}+\frac{\pi}4=\frac{(N-1)\pi}{4}>0$,
then there exist $n_0\in\N$,  $t_n\in[2^{n},2^{n}+2\pi]$ such that, letting $x_n=t_n e_1$
with $e_1=(1,0,\cdots,0)\in\R^N$,   for some $n\geq n_0$,
$$\frac34 2^{-\frac{N-1}{2}n}<\Psi(x_n):=\max_{t\in[2^{n},2^{n}+2\pi]} \Psi(te_1)<\frac54 2^{-\frac{N-1}{2}n}. $$
Note that
$$\lim_{n\to+\infty}(t_n-t_{n-1})2^{-n+1}=1.$$

Now we  denote  
\begin{equation}\label{5.1q}
 w_1(x)=\sum^{+\infty}_{n=n_0} a_n\Psi (x-x_n),\quad\forall\, x\in\R^N,
 \end{equation}
 where   $\{a_n\}_n$ is a sequence of nonnegative numbers $a_n=2^{\sigma n}$.

For $\sigma\in(-\frac{N-1}{2},\frac{N-1}{2})$,  we have that  $ \displaystyle\sum^{+\infty}_{n=n_0} a_n2^{-\frac{N-1}{2}n}<+\infty$ and
$$w_1(0)=\sum^{+\infty}_{n=n_0} a_n\Psi (x_n)<\frac54\sum^{+\infty}_{n=n_0} a_n2^{-\frac{N-1}{2}n}<+\infty.$$
Similarly, we get that $w_1(x)$ is bounded for any point $x$ in $\R^N$, so $w_1$ is well-defined in $\R^N$.

Note that
\begin{eqnarray*}
c|x_n|^{\sigma} \leq a_n\Psi(0)&<&w_1(x_n)
\\&<&a_n\Psi(0)+\sum^{n-1}_{m=n_0}a_{m}2^{-\frac{N-1}{2}(n-m)}+\sum^{+\infty}_{m=n+1}a_{m}2^{-\frac{N-1}{2}m}
\\& < &c|x_n|^{\sigma}+\sum^{n-1}_{m=n_0}a_{m}2^{-\frac{N-1}{2}(n-m)}+\sum^{+\infty}_{m=n+1}a_{m}2^{-\frac{N-1}{2}m},
\end{eqnarray*}
where
\begin{eqnarray*}
\sum^{n-1}_{m=n_0}a_{m}2^{-\frac{N-1}{2}(n-m)}   \leq  2^{-\frac{N-1}{2}n}\sum^{n-1}_{m=n_0} 2^{m(\frac{N-1}{2}+\sigma)}
  \leq  C 2^{\sigma n}=C|x_n|^{\sigma}
\end{eqnarray*}
and
\begin{eqnarray*}
\sum^{+\infty}_{m=n+1}a_{m}2^{-\frac{N-1}{2}m} & \leq & \sum^{+\infty}_{m=n+1} 2^{(\sigma-\frac{N+1}{2})m}     \\&\leq & c 2^{(\sigma-\frac{N+1}{2})n} = c |x_n|^{ \sigma-\frac{N+1}{2} }.
\end{eqnarray*}
Therefore, we have that
\begin{equation}\label{5.1q+1}
 \limsup_{|x|\to+\infty} |w_1(x)|(1+|x|)^{\sigma}\in(0,+\infty).
 \end{equation}

For $\sigma> \frac{N-1}{2}$,  let
\begin{equation}\label{5.1q+}
 w_2(x)=\sum^{+\infty}_{n=n_0} a_n\left[\Psi (x-x_n)-\Psi(x+x_n)\right],\quad\forall\, x\in\R^N,
 \end{equation}
 where   $\{a_n\}_n$ is a sequence of nonnegative numbers $a_n=2^{\sigma n}$.
Note that $w_2(0)=0$, so $w_2$ is well-defined in $\R^N$ and (\ref{5.1q+1}) holds.   
We complete the proof. \hfill$\Box$ \medskip

\noindent{\bf Proof of Proposition \ref{pr hom sol}.} Let $\psi$ be a solution of (\ref{1.2})
satisfying  that for some $\sigma>\frac{N-1}{2}$ and $c_1>0$,
$$ |\psi(x)|\leq c_1(1+|x|)^{-\sigma},\quad\forall\, x\in\R^N,$$
then
$$\lim_{r\to+\infty}\frac1r\int_{B_r(0)}|\psi(x)|^2 dx\leq \lim_{r\to+\infty}c_2r^{N-1-2\sigma }=0.$$
 From Rellich uniqueness, we have that
$\psi=0$. This means $\mathcal S_\sigma=\emptyset$ for $\sigma>\frac{N-1}{2}$.

For $\sigma=\frac{N-1}{2}$, $\mathcal S_\sigma$ contains the functions $\psi_k$ with $k\in \N$,
 defined by (\ref{1.1-2}).

For $-\infty<\sigma_1<\sigma_2\leq \frac{N-1}{2}$, it is easy to see that
$ {\mathcal S}_{\sigma_2}\subset{\mathcal S}_{\sigma_1}$. From Lemma \ref{lm 2.1-1}, we have that   for $-\infty<\sigma_2<\sigma_1< \frac{N-1}{2}$,
$$ {\mathcal S}_{\frac{N-1}{2}}\varsubsetneqq {\mathcal S}_{\sigma_1}\varsubsetneqq{\mathcal S}_{\sigma_2}.$$
We complete the proof. \hfill$\Box$

\setcounter{equation}{0}
 
\section{Classification}

 In this section, we build the connection between the singular solutions of (\ref{eq 1.1}) and the very weak solutions of (\ref{eq 1.1w}). \smallskip

  \begin{lemma}\label{lm 2.1}
Assume that  $N\ge2$,  $p>1$,  $Q:\R^N \to \R$ is a nonnegative H\"older continuous function such that $Q(0)>0$ and $u$ is a  classical solution of (\ref{eq 1.1}). 
 
 Then    $u\in   L^1_{loc}(\R^N)$ and $Q|u|^{p-1}u\in   L^1_{loc}(\R^N)$.
  
  \end{lemma}
 {\bf Proof. } From the assumptions that continuous $Q(0)>0$ and
$$\lim_{|x|\to0^+} u(x)=+\infty,$$
 there exists $r_0>0$ such that 
\begin{equation}\label{2.1}
Q>\frac{Q(0)}2,\ \  u>0\ \ {\rm in} \ \, B_{r_0}\setminus\{0\}.
\end{equation}  In fact,
 $u$ is a positive classical solution of
\begin{equation}\label{eq 2.1}
\arraycolsep=1pt\left\{
\begin{array}{lll}
 \displaystyle\ \ \, -\Delta u -u=Q|u|^{p-1}u\quad
 &{\rm in}\ \, B_{r_0}\setminus\{0\},\\[2mm]
 \phantom{  }
 \displaystyle   \lim_{|x|\to0^+}u(x)=+\infty.
\end{array}\right.
\end{equation}

  By contradiction, we assume that $f:=Q|u|^{p-1}u\not\in L^1(B_{r_0})$. Since $f$ is  positive and continuous in $B_{r_0}\setminus\{0\}$.
  
Let $\{r_n\}_n\subset (0,r_0)$ be a sequence of strictly decreasing positive numbers converging to zero and for any $r_n$, we have that
$$
 \lim_{r\to0^+} \int_{B_{r_n}(0)\setminus B_r(0)}f(x) dx   =+\infty,
$$
then there exists $R_n\in (0,r_n)$ such that
$$
  \int_{B_{r_n}(0)\setminus B_{R_n}(0)}f dx  =n,
$$
{\it Case of $\mu\ge0$.} Let $\delta_n=\frac1n \Gamma_\mu f\chi_{B_{r_n}(0)\setminus B_{R_n}(0)}$, then the problem
$$
 \arraycolsep=1pt\left\{
\begin{array}{lll}
 \displaystyle  -\Delta v = \delta_n\quad
   &{\rm in}\;\;  B_{r_0},\\[1mm]
 \phantom{ -\Delta    }
 \displaystyle  v= 0\quad  &{\rm   on}\;\; \partial B_{r_0}
 \end{array}\right.
$$
has a unique positive solution  $v_n$  satisfying (in the usual sense)
$$\int_{B_{r_0}} v_n  (-\Delta)\xi dx=\int_{B_{r_0}} \delta_n \xi dx,\quad\forall\, \xi\in C^{1.1}_0(B_{r_0}).$$
For any $\xi\in C^{1.1}_0(B_{r_0})$, we have that
$$\int_{B_{r_0}} w_n  (-\Delta) \xi\,  dx  =\int_{\Omega} \delta_n \xi\, dx\to \xi(0)\quad{\rm as}\;\; n\to+\infty. $$
Therefore,   for any compact set $\mathcal{K}\subset B_{r_0}\setminus \{0\}$
$$\norm{v_n-\frac1{c_N}G_{r_0}}_{C^1(\mathcal{K})}\to 0\quad{\rm as}\quad  {n\to+\infty}.$$

So we fix a point $x_0\in \Omega\setminus \{0\}$, let $d_0=\frac{\min\{|x_0|,\, r_0-|x_0|\}}{2}$ and $\mathcal{K}=  \overline{B_{d_0}(x_0)}$, then
there exists $n_0>0$ such that for $n\ge n_0$,
\begin{equation}\label{2.1-+-1}
 w_n\ge \frac1{2c_\mu}G_{r_0}\quad{\rm in}\;\; \mathcal{K},
\end{equation}
where $G_{r_0}$ is the weak solution of 
$$-\Delta u=\delta_0\quad {\rm in}\ \, B_{r_0},\qquad  u=0\quad {\rm on}\ \, \partial B_{r_0}.$$

Let $w_n$ be the solution   of
$$
 \arraycolsep=1pt\left\{
\begin{array}{lll}
 \displaystyle  -\Delta u = n\delta_n\quad
   &{\rm in}\ \,  B_{r_0},\\[1mm]
 \phantom{ -\Delta  }
 \displaystyle  u= 0\quad   &{\rm   on}\ \, \partial B_{r_0},
 \end{array}\right.
$$
then the comparison principle implies that $  w_n\ge nv_n\; {\rm in}\; B_{r_0}.$
Together with (\ref{2.1-+-1}), we derive that
$$v_n\ge  \frac n{2c_N}G_{r_0}\quad {\rm in}\;\;  \mathcal{K}.$$
Then by comparison principle again, we have that
$$u(x_0)\ge \frac n{2c_N}G_{r_0} \to+\infty\quad{\rm as}\;\; n\to+\infty,$$
which contradicts that $u$ is classical solution of (\ref{eq 1.1}).  

Therefore, we have that $Q|u|^{p-1}u\in L^1(B_{r_0})$ and $u\in L^1(B_{r_0})$
by the fact that $Q>\frac{Q(0)}{2}$ and $p>1$.
  \hfill$\Box$\medskip

The following estimates is important for our analysis of singularity and regularity.
\begin{proposition}\cite[Chapter V]{S} ( Proposition 2.2 in \cite{CQ} with $\alpha=1$) \label{embedding} 
Let $h\in L^s(\Omega)$ with $s\ge1$, then
there exists $c_{3}>0$ such that

\noindent$(i)$ if $\frac1s<\frac{2  }N$,
\begin{equation}\label{a 4.1}
\|\mathbb{G}_{r_0}[h]\|_{C^{\min\{1, 2 -N/s\}} (B_{r_0})}\le c_{3}\|h\|_{L^s(B_{r_0})};
\end{equation}

\noindent$(ii)$ if $\frac1s\le \frac1r+\frac{2}N$ and $s>1$,
\begin{equation}\label{a 4.2}
\|\mathbb{G}_{r_0}[h]\|_{L^r(B_{r_0})}\le c_{3}\|h\|_{L^s(B_{r_0})};
\end{equation}

\noindent$(iii)$ if $1<\frac1r+\frac{2 }N$,
\begin{equation}\label{a 4.02}
\|\mathbb{G}_{r_0}[h]\|_{L^r(B_{r_0})}\le c_{3}\|h\|_{L^1(B_{r_0})}.
\end{equation}
\end{proposition}

Now we are in the position to prove Theorem \ref{teo 1}. \medskip

\noindent{\bf Proof of Theorem \ref{teo 1}.}   Let $r_0>0$ satisfy (\ref{2.1})
 and then
 $u$ is a positive classical solution of
 $$
\arraycolsep=1pt\left\{
\begin{array}{lll}
 \displaystyle\ \, -\Delta u =u+ Q|u|^{p-1}u\quad
 &{\rm in}\quad B_{r_0}\setminus\{0\},\\[2mm]
 \phantom{  }
 \displaystyle   \lim_{|x|\to0^+}w(x)=+\infty.
\end{array}\right.
$$

Define the operator $L$  by the following
$$
L(\xi):=\int_{B_{r_0}} [u(-\Delta \xi-\xi) - Q|u|^{p-1}u\xi]\,dx,\quad \forall\xi\in C^\infty_c(B_{r_0}).
$$
First we claim that for any $\xi\in C^\infty_c(\R^N)$ with
the support in $B_{r_0}\setminus\{0\}$,
$$L(\xi)=0.$$
In fact, since $\xi\in C^\infty_c(\R^N)$ has the support in $B_{r_0}\setminus\{0\}$, then there exists $r\in(0,r_0/2)$ such that
$\xi=0$ in $B_r(0)$ and  then
\begin{eqnarray*}
L(\xi)
&=&  \int_{B_{r_0}\setminus B_r}[ u( -\Delta  \xi-\xi)  -Q|u|^{p-1}u\xi]\,dx
\\&  = &  \int_{B_{r_0}\setminus B_r}(-\Delta u-u  -Q|u|^{p-1}u )\xi\,dx
\\&=&0.
\end{eqnarray*}
Now let $\eta_0:\R^N\to[0,1]$ be a smooth radially symmetric function such that 
$\eta_0=1$ in $B_1$ and $\eta_0=0$ in $B_2^c$. 

By Lemma \ref{lm 2.1} we have that  $u\in L^1_{loc}(\R^N)$  and $Q|u|^{p-1}u\in L^1_{loc}(\R^N)$, then from Theorem 1 in \cite{BL},    it implies that
$$
 L=  k \delta_0\quad {\rm for\ some\ }\ k\ge0,
$$
 that is,
\begin{equation}\label{3.3-1}
L(\xi)= \int_{B_{r_0}} \left[u(-\Delta \xi-\xi) -Q|u|^{p-1}u\xi\right]\,dx=  k \xi(0),\quad  \quad \forall \,\xi\in C^{1.1}_c (B_{r_0}).
\end{equation}
Let $\Gamma$ be a Harmonic function such that $\Gamma=u$ on $\partial B_{r_0}$
then $\Gamma\in C^2(B_{r_0})$ and $u_0=u-\Gamma$ is a weak nonnegative solution of $$
\arraycolsep=1pt\left\{
\begin{array}{lll}
 \displaystyle -\Delta w =w+\Gamma+ Q|w+\Gamma|^{p-1}(w+\Gamma)+k\delta_0\quad
 &{\rm in}\quad B_{r_0} ,\\[2mm]
 \phantom{ -\Delta  }
 \displaystyle  w=0\quad &{\rm on}\quad \partial B_{r_0}
\end{array}\right.
$$for some $k\ge0$.

 When  $k=0$,  then
$$w= \mathbb{G}_{r_0}[f(w)], $$
where 
$$f(w)=w+\Gamma+ Q|w+\Gamma|^{p-1}(w+\Gamma).$$
 Here we have that for some $c_4>0$
$$|f(w)|\leq c_4 (w^p+1).$$

Let $N\geq 3$ and  we infer by Proposition \ref{embedding} that   $u^p\in L^{t_0}(B_{r_0})$ with $t_0=\frac{1}{2}[1+\frac1p\frac{N}{N-2}]>1$,   then use Proposition \ref{embedding} again
$u\in L^{t_1 p}(B_{r_0})$ and  $u^p\in L^{t_1}(B_{r_0})$ with
$$t_1=\frac1p\frac{N}{N-2t_0}t_0>t_0.$$
If $t_1>\frac1{2} Np$, by Proposition \ref{embedding},
$u\in L^{\infty}(B_{r_0})$ and then it could be improved that $u$ is a classical solution of $$
\arraycolsep=1pt\left\{
\begin{array}{lll}
 \displaystyle -\Delta w =w+\Gamma+ Q|w+\Gamma|^{p-1}(w+\Gamma)\quad
 &{\rm in}\quad B_{r_0} ,\\[2mm]
 \phantom{ -\Delta  }
 \displaystyle  w=0\quad &{\rm on}\quad \partial B_{r_0}
\end{array}\right.
$$

 If $t_1<\frac1{2} Np$, we proceed as above.
By Proposition \ref{embedding}, $u\in L^{t_2p}(B_{r_0})$, where
$$t_2=\frac1p\frac{Nt_1}{N-2  t_1 }>\frac1p\frac{N}{N-2 t_0 }t_1=\left(\frac1p\frac{N}{N-2  t_0}\right)^2t_0.$$
Inductively, let us define
$$t_m=\frac1p\frac{Nt_{m-1}}{N-2 t_{m-1}} >\left(\frac1p\frac{N}{N-2  t_0}\right)^m t_0\to+\infty\quad{\rm as}\quad m\to+\infty.$$
Then there exists $m_0\in\N$ such that
$$t_{m_0}>\frac{1}{2}Np$$
and by Proposition \ref{embedding} part $(i)$,
$$w\in L^\infty(B_{r_0}),$$
 then  $u\in L^\infty$, which contradicts    the assumption that $\displaystyle\lim_{|x|\to0^+} u(x)=+\infty.$  
 When $N=2$, we can choose $t_1> p$ and a contradiction is derived. 

Therefore,  $k>0$.\smallskip

{\it  Now we show the singularity of $u$ for $k>0$. }
Note that 
 \begin{equation}\label{13.2.1a}
w=\mathbb{G}_{r_0}[f(w)]+k\mathbb{G}_{r_0}[\delta_0],
\end{equation}
 letting
$$w_1= \mathbb{G}_{r_0}[f(w)]\quad {\rm and}\quad  \Upsilon_0=k\mathbb{G}_{r_0}[\delta_0],$$
then by  Young's inequality,
\begin{equation}\label{13.2.1b}
f(w)\le c_5\left(1+w_1^p+  \Upsilon_0^p\right).
\end{equation}
 By the definition of $w_1$ and (\ref{13.2.1b}), we obtain that
\begin{equation}\label{13.2.1c}
 w_1 \le c_5 \mathbb{G}_{r_0}[w_1^p]+ \Upsilon_1,
\end{equation}
where $w_1\in L^s(B_{r_0})$  for any  $s\in (1,  \frac{N}{N-2\alpha})$ and
$$ \Upsilon_1=c_5 \mathbb{G}_{r_0}[\Upsilon_0^p].$$
Denoting  $\mu_1=2+(2-N)p$, then for $0<|x|< \frac{r_0}2$,
$$\Upsilon_1(x)  \le
\left\{ \arraycolsep=1pt
\begin{array}{lll}
  c_6|x|^{\mu_1 } \quad
 &{\rm if}\quad \mu_1<0,\\[2mm]
 \displaystyle   -c_6\log |x| \quad
 &{\rm if}\quad \mu_1=0,\\[2mm]
 c_6 \quad
 &{\rm if}\quad \mu_1>0.
\end{array}
\right. $$
If $\mu_2\le 0$ ( $\mu_1>0$ for $N=2$),  denote
$$w_2= c_6 \mathbb{G}_{r_0}[ w_1^p]. $$
Then $w_2\in L^s(\Omega)$ with $s\in[1,\frac{N}{N-2 })$, $w_1\le w_2+\Upsilon_1$ and
$$w_2\le c_6\left(\mathbb{G}_\alpha[ w_2^p] + \mathbb{G}_\alpha[ \Upsilon_1^p]\right).$$
Let $\mu_2=\mu_1 p +2 $, then $\mu_2>\mu_1$ and for $0<|x|<\frac{r_0}2$,
$$\Upsilon_2(x):=c_6 \mathbb{G}_{r_0}[ \Upsilon_1^p](x)\le
\left\{ \arraycolsep=1pt
\begin{array}{lll}
  c_7|x|^{\mu_2 } \quad
 &{\rm if}\quad \mu_2<0,\\[2mm]
 \displaystyle   -c_7\log |x| \quad
 &{\rm if}\quad \mu_2=0,\\[2mm]
 c_7  \quad
 &{\rm if}\quad \mu_2>0.
\end{array}
\right. $$

Inductively, we assume that
$$w_{n-1}\le c_{n-1} \mathbb{G}_{r_0}[ w_{n-1}^p]+c_{n-1}  \mathbb{G}_\alpha[ \Upsilon_{n-2}^p], $$
where $c_{n-1}>0$, $w_{n-1}\in L^s(B_{r_0})$ for $s\in[1,\frac{N}{N-2\alpha})$, $\Gamma_{n-2}(x)\le |x|^{\mu_{n-2}}$ for  $\mu_{n-2}<0$.

Let
$$w_n= c_{n-1}  \mathbb{G}_{r_0}[ w_{n-1}^p],\quad \quad \Upsilon_{n-1}=c_{n-1} \mathbb{G}_{r_0}[ \Upsilon_{n-2}^p] $$
and
 $$\mu_{n-1}=\mu_{n-2} p +2 ,$$
then $u_n\in L^s(B_{r_0})$ for $s\in[1,\frac{N}{N-2})$ and  for $0<|x|<\frac{r_0}2$ and $c_n>0$,
$$\Upsilon_{n-1}(x):= \mathbb{G}_{r_0}[  \Upsilon_{n-2}^p](x)\le
\left\{ \arraycolsep=1pt
\begin{array}{lll}
  c_n|x|^{\mu_{n-1} } \quad
 &{\rm if}\quad \mu_{n-1}<0,\\[2mm]
 \displaystyle   -c_n\log |x| \quad
 &{\rm if}\quad \mu_{n-1}=0,\\[2mm]
 c_n  \quad
 &{\rm if}\quad \mu_{n-1}>0.
\end{array}
\right. $$
We observe that
\begin{eqnarray*}
 \mu_{n-1}-\mu_{n-2}= p(\mu_{n-2}-\mu_{n-3})&=&p^{n-3}(\mu_2-\mu_1)  \\
   &\to&+\infty\quad{\rm as}\quad n\to+\infty.
\end{eqnarray*}
Then there exists $n_2\ge 1$ such that
 $$ \nu_{n_2-1}>0\quad{\rm and}\quad \nu_{n_2-2}\le0$$
and
\begin{equation}\label{3.2.1d1-}
w\le w_{n_2}+\sum^{n_2-1}_{i=1} \Upsilon_i +\Upsilon_0,
\end{equation}
where $\Upsilon_i\le c|x|^{\mu_i}$ and
$$w_{n_2}\le c_{n_2} (\mathbb{G}_{r_0}[ w_{n_2}^p]+1).$$

We next claim that $w_{n_2}\in L^\infty(\Omega)$. Since $w_{n_2}\in L^s(\Omega)$ for  $s\in[1,\frac{N}{N-2 })$,
letting
$$t_0=\frac12(1+\frac1p\frac{N}{N-2 })\in (1,\frac{N}{N-2}),$$
 then $\frac1p\frac{N}{N-2  t_0 }>1$ and
by Proposition \ref{embedding}, we have that
$w_{n_2}\in L^{t_1}(B_{r_0})$ with
$$t_1=\frac1p\frac{Nt_0}{N-2  t_0 }.$$
Inductively, it implies by $w_{n_2}\in L^{t_{n-1}}(B_{r_0})$ that
$w_{n_2}\in L^{t_{n}}(B_{r_0})$
with
$$t_n=\frac1p\frac{Nt_{n-1}}{N-2  t_{n-1}} >\left(\frac1p\frac{N}{N-2  t_0}\right)^n t_0\to+\infty\quad{\rm as}\quad n\to\infty.$$
Then there exists $n_3\in\N$ such that
$$s_{n_3}>\frac{Np}{2 }$$
and by Proposition \ref{embedding} part $(i)$, it infers that
$$w_{n_2}\in L^\infty(B_{r_0}).$$

Therefore, it implies by $w\ge \Upsilon_0 $ and (\ref{3.2.1d1-}) that
$$\lim_{|x|\to0^+} w(x)|x|^{N-2}=c_Nk\ \ {\rm if}\ \ N\geq 3,\qquad \lim_{|x|\to0^+} \frac{w(x)}{-\ln |x|} =c_Nk\ \ {\rm if}\ \ N=2,$$
where $c_N>0$ is the normalized constant. 

For $C^\infty_c (\R^N)$, we can divided $\xi$ into $\xi_1+\xi_2$, where $\xi_1,\xi_2$ are smooth function such that  
$${\rm supp}(\xi_1)\subset B_{r_0}\quad {\rm supp}(\xi_2)\subset B_{\frac{r_0}2}^c.$$
Since $u$ is a classical solution of (\ref{eq 1.1}), then we have that 
$$
\mathcal{L} (\xi_2)= \int_{\R^N} \left[u(-\Delta) \xi_2-u\xi_2 -Q|u|^{p-1}u\xi_2\right]\,dx= 0.
$$
which, together with (\ref{3.3-1})  replaced $\xi$ by $\xi_1$, implies that 
$$
 \int_{\R^N} \left[u(-\Delta) \xi-u\xi - Q|u|^{p-1}u\xi\right]\,dx=  k \xi(0),\quad  \forall\, \xi\in C^\infty_c (\R^N).
$$
Since any $\xi\in C^{1,1}_c(\R^N)$ could be  approximated by a sequence of functions in $C^\infty_c (\R^N)$, so (\ref{3.3-1}) holds for any $\xi\in C^{1,1}_c(\R^N)$.
\hfill$\Box$\medskip



\setcounter{equation}{0}
\section{Existence of weak solutions}

We first provide some important estimates from the convolution  by fundamental solution of Helmholtz operator.

   \begin{lemma}\label{lm 4.1}
For 
  $$N\geq 2,\ \  \tau>\frac{N+1}2,\ \  
 \theta\in \big(-1,\, N-2\big)\setminus\{0\},$$
 let  $U_{\theta,\tau}:\R^N\setminus\{0\}\to\R$ be a  function in $L^1_{loc}(\R^N)$ such that for some $c_8>0$
\begin{equation}\label{3.4-0--}
  |U_{\theta,\tau} (x)|\leq c_8|x|^{-\theta-2}(1+|x|)^{-\tau+\theta+2},\quad {\rm a.e.\ in}\ \ \R^N\setminus\{0\}.
   \end{equation}
  
Then there exists $c_9>0$ such that 
\begin{equation}\label{3.4}
 | (\Phi\ast U_{\theta,\tau})(x) |  \le c_9 (1+ |x|^{-\theta})(1+|x|)^{\max\left\{-\frac{N-1}{2},\frac{1+N}2-\tau\right\}+\theta},\quad \forall\, x\in\R^N
\end{equation}
and
\begin{equation}\label{3.5}
 | (\nabla\Phi\ast U_{\theta,\tau})(x) |  \le c_9  |x|^{-\theta-1}(1+|x|)^{\max\left\{-\frac{N-1}{2},\frac{1+N}2-\tau\right\}+\theta+1},\quad \forall\, x\in\R^N.
\end{equation}

\end{lemma}
{\bf Proof.} In the following, the letter $c_j>0$ always denotes constants which only depends on $N$, $\alpha$ and $k$. We observe that 
\[
|\Phi(x)|\leq
\left\{\begin{array}{lll}
c_{10}\, |x|^{2-N}\quad\ \ &{\rm if }\ N\geq 3,\\[1.5mm]
 \phantom{   }
c_{10}\, \log \frac2{|x|}   & {\rm if }\ N=2,\,
\end{array}\right.
\qquad |\nabla\Phi|\leq c_{11}|x|^{1-N}
\quad\, {\rm for} \ \ 0<|x|\leq 1 \]
and
\[ |\Phi(x)|,\, |\nabla\Phi|\leq c_{12}\, |x|^{\frac{1-N}2} \quad  {\rm if }\ \  |x|>1. \]

Note that  $$ U_{\theta,\tau}(x)   \le  2|x|^{-2-\theta}\chi_{B_1}(x) + 2|x|^{-\tau}\chi_{B_1^c}(x),\quad x\in \R^N\setminus \{0\},$$
where 
$B_r(z)$ is the ball centered at $z$ with radius $r>0$, $B_r(0)=B_r$, $A^c=\R^N\setminus A$ and
$\chi_{A}$ is characteristic function of $A$, i.e. $\chi_A(x)=1$ for $x\in A$ and
$\chi_{A^c}(x)=0$ for $x\in A^c$.

Let
 $$U_1(x)=|x|^{-2-\theta }\chi_{B_1}(x),\quad U_2(x)= |x|^{-\tau}\chi_{B_1^c}(x),\quad\forall\, x\in\R^N\setminus\{0\}.$$
 From the fact that $-\Delta (|x|^{-\theta})=\theta (N-2-\theta )|x|^{-2-\theta }$ in $\mathbb{R}^N\setminus\{0\}$, we can deduce
$$C_N\int_{\R^N}\frac{|y|^{-\theta -2} }{|x-y|^{N-2}}dy=\frac{1}{\theta (N-2-\theta )}|x|^{-\theta}$$
and
$$C_2\int_{\R^2} |y|^{-\theta -2}  (-\ln |x-y|) dy\leq c_{13}.$$

For $x\in B_2\setminus\{0\}$, we have that for $N\geq 3$ and $\theta\in(0,N-2)$
 \begin{eqnarray*}
 |(\Phi\ast U_1)(x)| \leq    C_N\int_{B_1}\frac{|y|^{-\theta -2} }{|x-y|^{N-2}}dz\leq  \frac{1}{\theta (N-2-\theta )}|x|^{-\theta },
\end{eqnarray*}
and for $N=2$ and $\theta\in(-1,\, 0)$
\begin{eqnarray*}
 |(\Phi\ast U_1)(x)|&\leq &  C_N\int_{B_1} |y|^{-\theta -2} (-\ln|x-y|)dz\leq c_{14}.
\end{eqnarray*}

When $x\in \R^N\setminus B_2(0)$,
$$
  |\Phi\ast U_1)(x) |\leq c (|x|-1)^{\frac{1-N}{2}} \int_{B_1} |y|^{-\theta -2}dy
\leq c_{15}  |x|^{\frac{1-N}{2}}.
$$
 
Thus
\begin{equation}\label{3.2}
 \big|(\Phi\ast U_1)(x) \big| \le  c_{16}  |x|^{-\theta }(1+|x|)^{\frac{1-N}2+\theta},\quad  \forall\, x\in\R^N\quad{\rm for}\ \, N\geq 3
  \end{equation}
  and
 \begin{equation}\label{3.2 N=2}
 \big|(\Phi\ast U_1)(x) \big| \le  c_{17}  (1+|x|)^{-\frac{1}2},\quad  \forall\, x\in\R^2.
  \end{equation} 
  
   \smallskip

  Note that for $x\in B_4\setminus\{0\}$, we have that for $N\geq 3$ 
  \begin{eqnarray*}
  |(\Phi\ast U_2)(x) |  &\leq & C_N\int_{B_8\setminus B_1}\frac{  |y|^{-\tau} }{|x-y|^{N-2} }dy + C_N\int_{  B_8^c }\frac{|y|^{-\tau} }{|x-y|^{\frac{N-1}2}}dy
  \\&\leq&     c_{18} \Big(1+\int_{ B_8^c}|y|^{-\tau-\frac{N-1}{2}}dz\Big)
\end{eqnarray*}
and for $N=2$, 
\begin{eqnarray*}
  |(\Phi\ast U_2)(x) |  \leq  C_N\int_{B_8\setminus B_1}  |y|^{-\tau} (-\ln |x-y|)_+ dy + C_N\int_{  B_8^c }\frac{|y|^{-\tau} }{|x-y|^{\frac{1}2}}dy
 \leq   c_{18},
\end{eqnarray*}
thus,
\begin{eqnarray*}
  |(\Phi\ast U_2)(x) |   \leq c_{19}.
\end{eqnarray*}
 
 For $x\in B_4^c$, there holds for $N\geq 3$,
\begin{eqnarray*}
  |(\Phi\ast U_2)(x) |   &\leq & \int_{B_{1}(x)}\frac{|y|^{-\tau} }{|x-y|^{N-2}}dy +\int_{B_{|x|/2}(x)\setminus B_1(x)}\frac{|y|^{-\tau} }{|x-y|^{\frac{N-1}{2}}}dy
  \\&& +\int_{B_{|x|/2}(0)\setminus B_1(0)}\frac{|y|^{-\tau} }{|x-y|^{\frac{N-1}{2}}}dy+ \int_{ \big( B_{|x|/2}(0) \cup B_{|x|/2}(x)\big)^c} |y|^{-\tau-\frac{N-1}{2}}dy\\
&\leq&c_{20}\big(|x|^{-\tau} +|x|^{-\tau+\frac{N+1}{2}}\big)    
\end{eqnarray*}
for $N= 2$,
\begin{eqnarray*}
  |(\Phi\ast U_2)(x) |   &\leq & \int_{B_{1}(x)} |y|^{-\tau} \ln|x-y| dy +\int_{B_{|x|/2}(x)\setminus B_1(x)}\frac{|y|^{-\tau} }{|x-y|^{\frac{1}{2}}}dy
  \\&& +\int_{B_{|x|/2}\setminus B_1}\frac{|y|^{-\tau} }{|x-y|^{\frac{1}{2}}}dy+ \int_{ ( B_{|x|/2} \cup B_{|x|/2}(x))^c} |y|^{-\tau-\frac{1}{2}}dy\\
&\leq&c_{21}(|x|^{-\tau} +|x|^{-\tau+\frac{3}{2}}),   
\end{eqnarray*}
where we used the fact  that $\tau>\frac{3}2$.

 Thus there exists $c_{22}>0$ such that 
 \begin{equation}\label{3.3}
  \big|(\Phi\ast U_2)(x) \big|  \le c_{22}  (1+|x|)^{\frac{N+1}2-\tau},\quad \forall\, x\in\R^N.
\end{equation}
Therefore,  (\ref{3.4})  follows  by (\ref{3.2})  and (\ref{3.3}) for $N\geq 3$,
and    by (\ref{3.2 N=2})  and (\ref{3.3}) for $N=2$. 

Similarly, we can have the following  estimates:
For $N\geq 2$, we have that 
\begin{equation}\label{3.2-nabla}
\big|(\nabla\Phi\ast U_1)(x) \big| \le  c_{23}  |x|^{-\theta-1 }(1+|x|)^{\frac{1-N}2+\theta+1},\quad  \forall\, x\in\R^N 
  \end{equation}
 and
  \begin{equation}\label{3.3-nabla}
  \big|(\nabla\Phi\ast U_2)(x)\big|   \le c_{24}  (1+|x|)^{\frac{1+N}2-\tau},\quad \forall\, x\in\R^N,
\end{equation}
which imply (\ref{3.5}). 
\hfill$\Box$

\smallskip

 Now we are in the position to show the existence of singular solution of (\ref{eq 1.1}).
 
 \medskip
 
 \noindent{\bf Proof of Theorem \ref{teo 2}. }  
   Denote
$$w_{\sigma}(x):=\Phi(x)+\phi_\sigma(x),\quad x\in\R^N\setminus\{0\},$$
where $\phi_\sigma\in \mathcal S_{\sigma}$ and 
$$\sigma\in\left[ \frac{1}{p-1}\Big(\frac{N+1}{2}-\alpha\Big),\, \frac{N-1}{2}\right],$$
which is a nonempty interval  by the assumption of $p$.

Then 
$w_{\sigma}$ is a solution weak solution of 
$$-\Delta w_{\sigma}-w_{\sigma}= \delta_0\quad{\rm in}\quad \R^N.$$
 Note that there is function $\psi_{\sigma}\in \mathcal{S}_{\sigma}^o$ such that  $w_{\sigma}=\Phi+\psi_\sigma$ is smooth in $\R^N\setminus\{0\}$ and  
$$\lim_{|x|\to0}w_{\sigma}(x)\cdot|x|^{N-2}=C_N\quad{\rm and}\quad 
 \limsup_{|x|\to+\infty}|w_{\sigma}(x)| |x|^{\sigma}= c_{25}. $$

To obtain  a solution of (\ref{eq 1.1w}), we turn to obtain the weak solution of 
\begin{equation}\label{eq 4.1}
    -\Delta v-v=Q |kw_\sigma+v|^{p-1}(kw_\sigma+v)  \quad
  {\rm in}\quad \R^N
\end{equation}
by considering the equivalent equation 
\begin{equation}\label{eq 4.2}
v=\Phi\ast \Big(Q |kw_\sigma+v|^{p-1}(kw_\sigma+v) \Big). 
\end{equation}

For $p<p^*_{_N}$, we have that  
$(2-N)p+2>2-N$, and now we fix
  $$\theta_p=
  \left\{\begin{array}{lll}
 \frac{2-N}2+\frac{(2-N)p+2}{2}\ \ &{\rm if }\ \ 2-(N-2)p\leq 0,\\[2mm]
 \phantom{   }
0\quad\ \ &{\rm if }\  \  2-(N-2)p>0 
\end{array}\right. 
$$ 
and   denote  
\[
W_p (x) = 
|x|^{\theta_p}(1+|x|)^{-  \sigma-\theta_p}  \quad {\rm for}\ \, x\in\R^N\setminus\{0\}.
\]

It is worth noting that 
\[W_p (x) =
 (1+|x|)^{- \sigma }\quad   {\rm if }\ \ 2-(N-2)p> 0.
 \]

 Now let's denote   
\begin{eqnarray*}
\mathcal{D}_{p,k}  = \Big\{v\in L^1(\R^N):  \   |v|\le   k W_p \ \ {\rm a.e.\ in}\ \, \R^N\Big\}
\end{eqnarray*}
and
$$\mathcal{T}v=  \Phi\ast  [Q |v+kw_\sigma|^{p-1}(v+kw_\sigma)],\quad\forall\, v\in \mathcal{D}_{p,k}.$$

Observe that for $N\geq 3$
\[Q |v+kw_\sigma |^{p} \leq   \left\{\begin{array}{lll}
c_{26}k^p|x|^{(2-N)p} \ \  &{\rm for }\ \ 0<|x|\leq 1\\[2mm]
 \phantom{   }
c_{26}k^p |x|^{-\alpha- \sigma p}\  \ &{\rm for }\ \ |x|> 1
\end{array}\right.\]
and
for $N\geq 3$
\[Q |v+kw_\sigma |^{p} \leq   \left\{\begin{array}{lll}
c_{26}k^p(1-\ln |x|)^p  \ \  &{\rm for }\ \ 0<|x|\leq 1\\[2mm]
 \phantom{   }
c_{26}k^p |x|^{-\alpha- \sigma p}\  \ &{\rm for }\ \ |x|> 1,
\end{array}\right.\]
then we apply  Lemma \ref{lm 4.1} to obtain  that
\begin{equation}\label{eq 4.1-1-1}
\Big|\Phi\ast  [Q |v+kw_\sigma|^{p-1}(v+kw_\sigma)]\Big|  \leq 
 \left\{\begin{array}{lll}
c_{27}k^p|x|^{\theta_p} \ \ &{\rm for }\ \ 0<|x|\leq 1,\\[2mm]
 \phantom{   }
c_{27}k^p |x|^{-\sigma_p}\ \   &{\rm for }\ \ |x|> 1,
\end{array}\right. 
\end{equation}
where  $c'>0$ is independent of $k$, 
$$\sigma_p=\min\Big\{\frac{N-1}{2},\alpha+ \sigma p-\frac{1+N}2\Big\}.$$ 

In order to show $\mathcal{T} \mathcal{D}_{p,k} \subset  \mathcal{D}_{p,k}$,
we need $$\sigma\leq   \sigma_p,$$ that is,
\begin{equation}\label{eq 4.1}
 \sigma\leq \alpha+\sigma p-\frac{1+N}2.
 \end{equation}
Note that   (\ref{eq 4.1}) holds by our choice of $\sigma$ i.e.
$$\sigma\in\Big[ \frac{1}{p-1}\Big(\frac{N+1}{2} -\alpha\Big) ,\,  \frac{N-1}{2}\Big]$$
and
$$\frac{1}{p-1}\Big(\frac{N+1}{2}-\alpha\Big)\leq \frac{N-1}{2},$$
 which is guaranteed by $p\geq p^\#_\alpha$.

By the assumption $p>1$, there exists $k^*>0$ such that for $k\in(0,k^*)$
$$
  \mathcal{T} \mathcal{D}_{p,k} \subset  \mathcal{D}_{p,k}.
$$
Note that for $v\in \mathcal{D}_{p,k}$,  
applying Lemma \ref{lm 4.1} again we  that
\[\Big|\nabla\Phi\ast  [Q |v+kw_\sigma|^{p-1}(v+kw_\sigma)]\Big|  \leq 
 \left\{\begin{array}{lll}
c_{28}k^p|x|^{\theta_p-1} \ \  &{\rm for }\ \ 0<|x|\leq 1,\\[2mm]
 \phantom{   }
c_{28}k^p |x|^{ -\sigma_p }\ \ &{\rm for }\ \ |x|> 1.
\end{array}\right.\]

So we have that $\mathcal{T} v\in W^{1.\infty}_{loc}(\R^N\setminus\{0\})$,
by the embeddings $W^{1,\infty}(B_{R}\setminus B_{\frac1R})\hookrightarrow   L^1( B_{R}\setminus B_{\frac1R})$   are  compact and
together with upper bound in $\mathcal{D}_{p,k}$, we obtain
  that $\mathcal{T}$ is a compact operator in $\mathcal{D}_{p,k}$.

Observing that $\mathcal{D}_{p,k}$ is a closed and convex set in $L^1(\R^N)$, we may apply Schauder
fixed point theorem to derive that there exists $v_k\in \mathcal{D}_{p,k}$ such that
$$\mathcal{T}v_k=v_k. $$
Since $ |v_k|\le  k W_p$, so $v_k$ is locally bounded in $\R^N\setminus\{0\}$,
then $u_k:=v_k+kw_\sigma$ satisfies 
$$\limsup_{|x|\to+\infty} |u(x)||x|^{\sigma}<+\infty.$$

Moreover, if $ \sigma_p>    \sigma,$  then (\ref{eq 4.1-1-1}) implies that 
$$|v_k(x)|=|\Phi\ast [Q |v_k+kw_\sigma|^{p-1}(v_k+kw_\sigma)]|\leq c_{29}k^p|x|^{-\sigma_p}\quad {\rm for}\ \, |x|>1$$
and  then
$$\limsup_{|x|\to+\infty} |u_k(x)-kw_\sigma(x)||x|^{\sigma_p}\in[0,+\infty).$$
By the standard interior regularity results, $u_k$ is a positive classical solution of (\ref{eq 1.1}).
From Theorem \ref{teo 1} it implies that $u_k$ is a distributional solution of (\ref{eq 1.1}). \hfill$\Box$

 \setcounter{equation}{0}
\section{Complex valued solutions}

 Let ${\bf u}$ be  a complex value solutions of (\ref{eq 1.1}), then  it could be written as 
$${\bf u}=u_1+{\rm i} u_2,$$
where $u_1,u_2$ are the isolated singular solution of system in real valued framework
\begin{equation}\label{eq 1.1-com}
\left\{
 \arraycolsep=1pt
\begin{array}{lll}
 \  -\Delta u_1-u_1=Q|{\bf u}|^{p-1} u_1 \quad   &\rm{in}\ \, \R^N\setminus\{0\},\\[2mm]
 \  -\Delta u_2-u_2=Q|{\bf u}|^{p-1} u_2 \quad   &\rm{in}\ \, \R^N\setminus\{0\},\\[2mm]
\displaystyle
\lim_{|x|\to0^+}u_1(x)=+\infty,
\end{array}\right.
 \end{equation}
 where $|{\bf u}|=(u_1^2+u_2^2)^{\frac{1}2}$.
 
Moreover, the complex valued weak solution of (\ref{eq 1.1w}) could be written 
\begin{equation}\label{eq 1.1 comw}
\left\{
 \arraycolsep=1pt
\begin{array}{lll}
\displaystyle  -\Delta u_1-u_1=Q|{\bf u}|^{p-1}u_1+k\delta_0 \quad   &\rm{in}\ \, \mathcal{D}'(\R^N) ,\\[2mm]
 \displaystyle  -\Delta u_2-u_2=Q|{\bf u}|^{p-1}u_2 \quad   &\rm{in}\ \, \mathcal{D}'(\R^N).
\end{array}\right.
 \end{equation}
It is worthing noting that, in the complex valued framework,  the isolated singular solution of (\ref{eq 1.1}) is no longer to be classified   by the Dirac mass in (\ref{eq 1.1}) in the Serrin's subcritical case.  In fact,  it should be classified by
\begin{equation}\label{eq 1.1w-com}
  -\Delta   {\bf u}-{\bf u} =Q|{\bf u}|^{p-1}{\bf u}+k_1\delta_0+k_2 {\rm i} \delta_0\ \  \, \rm{in}\ \, \R^N
 \end{equation}
 or in the form of system (in the real valued framework)
$$
\left\{
 \arraycolsep=1pt
\begin{array}{lll}
\displaystyle  -\Delta u_1-u_1=Q|{\bf u}|^{p-1}u_1+k_1\delta_0 \quad   &\rm{in}\ \, \mathcal{D}'(\R^N) ,\\[2mm]
 \displaystyle  -\Delta u_2-u_2=Q|{\bf u}|^{p-1} u_2+k_2 \delta_0  \quad   &\rm{in}\ \, \mathcal{D}'(\R^N).
\end{array}\right.
$$

However, the weak solutions of   (\ref{eq 1.1w})  are classical isolated singular solutions of  (\ref{eq 1.1}). Therefore,  our interest in this section is to obtain weak solutionx  (\ref{eq 1.1w}) in the Serrin's subcritical case and the result states as follows.

  \begin{theorem}\label{teo 3}
Assume that $N\geq  2$,  the potential $Q$  verifies $( \mathbb{Q}_\alpha)$
 with $\alpha$ satisfying  (\ref{req 3}). 
 Let $p^\#_\alpha,\, p^*_{_N}$ be defined in (\ref{serrin 1}) and (\ref{serrin 2}) respectively.
 
 Then for 
 $$p\in \big(1,\,+\infty\big)\cap \big[p^\#_\alpha,\, p^*_{_N} \big)$$
and   
$\psi_\sigma\in \mathcal{S}_\sigma$ 
with 
$$\sigma\in\Big[ \frac{1}{p-1}\Big(\frac{N+1}{2} -\alpha\Big) ,\,  \frac{N-1}{2}\Big],$$
 there exists $k^*>0$
  such that  for $ k\in(0, k^*)$, problem (\ref{eq 1.1w}) admits a complex valued solution ${\bf u}_{k}$
  in a integral form
  $${\bf u}_k=\Phi_c\ast (Q|{\bf u}_k|^{p-1}{\bf u}_k)+k(\Phi_c+\psi_\sigma)\quad{\rm in}\ \ \R^N. $$
 
  \end{theorem}

 In order to obtain the complex valued solution, we need the following estimate.

  \begin{lemma}\label{cr 4.1}
Let 
  $$N\geq 2,\ \  \tau>\frac{N+1}2,\ \  
 \theta\in \big(-1,\, N-2\big)\setminus\{0\},$$
then for   ${\bf U}_{\theta,\tau}\in L^1(\R^N,\C)$  verifying  
  $$  |{\bf U}_{\theta,\tau} (x)|\leq |x|^{-\theta-2}(1+|x|)^{-\tau+\theta+2},\quad{\rm a.e.\ for} x\in\R^N\setminus\{0\},$$
 there holds for $c_{30}>0$   
\begin{equation}\label{3.4-c}
 | (\Phi_c\ast {\bf U}_{\theta,\tau})(x) |  \le c_{30} (1+ |x|^{-\theta})(1+|x|)^{\max\left\{-\frac{N-1}{2},\frac{1+N}2-\tau\right\}+\theta},\quad \forall\, x\in\R^N
\end{equation}
and
\begin{equation}\label{3.5-c}
 | (\nabla\Phi_c\ast {\bf U}_{\theta,\tau})(x) |  \le c_{30}  |x|^{-\theta-1}(1+|x|)^{\max\left\{-\frac{N-1}{2},\frac{1+N}2-\tau\right\}+\theta+1},\quad \forall\, x\in\R^N.
\end{equation}

\end{lemma}
{\bf Proof. } Note that $\Phi_c=\Phi+{\rm i}\Psi$, ${\bf U}_{\theta,\tau}=U _{\theta,\tau,1}+{\rm i} U_{\theta,\tau,2}$
 $$\Phi\ast {\bf U}_{\theta,\tau}=\Phi\ast    U _{\theta,\tau,1}+{\rm i}\Phi\ast   U_{\theta,\tau,2}$$ and 
 $$\nabla \Phi\ast {\bf U}_{\theta,\tau}=\nabla\Phi\ast   U _{\theta,\tau,1}+{\rm i}\nabla \Phi\ast   U _{\theta,\tau,2},$$ 
then  the related estimates  could see Lemma \ref{lm 4.1}.  So we only need to prove 
that  there exists $c_{32}>0$  such that for $j=1,2$
\begin{equation}\label{3.4-2}
 | (\Psi\ast {\bf U}_{\theta,\tau,j})(x) |  \le c_{32}  (1+|x|)^{\frac{N+1}{2}-\tau},\quad \forall\, x\in\R^N
\end{equation}
and
\begin{equation}\label{3.5-2}
 | (\nabla\Psi\ast {\bf U}_{\theta,\tau,j})(x) |  \le c_{32}  (1+|x|)^{\frac{N+1}{2}-\tau },\quad \forall\, x\in\R^N.
\end{equation}
 
Indeed, note that  $\Psi$ and $\nabla\Psi$ are smooth in $\R^N$ and have the asymptotic behavior at infinity controlled by $(1+|x|)^{-\frac{N-1}{2}}$.
Let
 $$U_1(x)=|x|^{-2-\theta }\chi_{B_1}(x),\quad U_2(x)= |x|^{-\tau}\chi_{B_1^c}(x),\quad\forall\, x\in\R^N\setminus\{0\}.$$
 Similar computations in Lemma \ref{lm 4.1}, we have that 
\begin{equation}\label{3.2--2}
 \big|(\Psi\ast U_1)(x) \big| \le  c_{32}  (1+|x|)^{\frac{1-N}2},\quad  \forall\, x\in\R^N 
  \end{equation}
  and
   \begin{equation}\label{3.2--2-}
 \big|(\Psi\ast U_2)(x) \big| \le  c_{32}  (1+|x|)^{\frac{1+N}2-\tau},\quad  \forall\, x\in\R^N,  
  \end{equation}
which imply (\ref{3.4-2}). 

Note that  (\ref{3.2--2}) and (\ref{3.2--2-}) hold true replacing $\Psi$ by $\nabla \Psi$, then (\ref{3.5-2}) holds true. \hfill$\Box$\medskip

\noindent{\bf Proof of Theorem \ref{teo 3}. }  
Set
$${\bf w}_{\sigma}(x):=\Phi_c(x)+\phi_\sigma(x),\quad x\in\R^N\setminus\{0\},$$
where $\phi_\sigma\in \mathcal S_{\sigma}$ and 
$$\sigma\in\left[  \frac{1}{p-1}\Big(\frac{N+1}{2}-\alpha\Big) ,\, \frac{N-1}{2}\right],$$
It is  a complex valued solution weak solution of 
$$-\Delta {\bf w}_{\sigma}-{\bf w}_{\sigma}= \delta_0\quad{\rm in}\ \, \R^N.$$

To obtain  a complex valued  solution of (\ref{eq 1.1w}),  we would consider the complex valued solution of 
\begin{equation}\label{eq 4.2-c}
{\bf v}=\Phi_c\ast \Big(Q |k{\bf w}_\sigma+{\bf v}|^{p-1}(k{\bf w}_\sigma+{\bf v}) \Big). 
\end{equation}
Recall that  
  $$\theta_p=
  \left\{\begin{array}{lll}
 \frac{2-N}2+\frac{(2-N)p+2}{2}\ \ &{\rm if }\ \ 2-(N-2)p\leq 0,\\[2mm]
 \phantom{   }
0\quad\ \ &{\rm if }\  \  2-(N-2)p>0 
\end{array}\right. 
$$ 
and   
\[
W_p (x) = 
|x|^{-\theta_p}(1+|x|)^{-  \sigma+\theta_p}  \quad {\rm for}\ \, x\in\R^N\setminus\{0\}.
\]

Re-denote   
\begin{eqnarray*}
 \mathcal{D}_{p,k}  = \Big\{{\bf v}\in L^1(\R^N;\C):  \   |{\bf v}|\le   k W_p \ \ {\rm a.e.\ in}\ \, \R^N\Big\}
\end{eqnarray*}
and
$$\mathcal{T}_c{\bf v}=  \Phi_c\ast  [Q |{\bf v}+k{\bf w}_\sigma|^{p-1}({\bf v}  +k{\bf w}_\sigma)],\quad\forall\, {\bf v}\in \mathcal{D}_{p,k}.$$

Observe that 
\[|Q| |{\bf v}+k{\bf w}_\sigma |^{p} \leq   \left\{\begin{array}{lll}
c_{33}k^p|x|^{(2-N)p} \ \  &{\rm for }\ \ 0<|x|\leq 1\\[2mm]
 \phantom{   }
c_{33}k^p |x|^{-\alpha- \sigma p}\  \ &{\rm for }\ \ |x|> 1,
\end{array}\right.\]
then we apply  Lemma \ref{cr 4.1} to obtain  that
\begin{equation}\label{eq 4.1-1-1-c}
\Big|\Phi_c\ast  [Q |{\bf v}+k{\bf w}_\sigma|^{p-1}({\bf v}+k{\bf w}_\sigma)]\Big|  \leq 
 \left\{\begin{array}{lll}
c_{34}k^p|x|^{\theta_p} \ \ &{\rm for }\ \ 0<|x|\leq 1,\\[2mm]
 \phantom{   }
c_{34}k^p |x|^{-\sigma_p}\ \   &{\rm for }\ \ |x|> 1,
\end{array}\right. 
\end{equation}
where  $c_{34}>0$ is independent of $k$, 
$$\sigma_p=\min\Big\{\frac{N-1}{2},\alpha+ \sigma p-\frac{1+N}2\Big\}.$$ 

In order to show $\mathcal{T}_c \mathcal{D}_{p,k} \subset  \mathcal{D}_{p,k}$,
we need $$\sigma \leq   \sigma_p,$$ that is,
\begin{equation}\label{eq 4.1-c}
 \sigma\leq \alpha+\sigma p-\frac{1+N}2 .
 \end{equation}
Note that the first inequality of (\ref{eq 4.1}) holds by our choice of $\sigma$ i.e.
$$\sigma\in\Big[ \frac{1}{p-1}\Big(\frac{N+1}{2} -\alpha\Big) ,\,  \frac{N-1}{2}\Big]$$
by $p\geq p^\#_\alpha$.

By the assumption $p>1$, there exists $k^*>0$ such that for $k\in(0,k^*)$
$$
  \mathcal{T}_c \mathcal{D}_{p,k} \subset  \mathcal{D}_{p,k}.
$$
Note that for $v\in \mathcal{D}_{p,k}$,  
applying Lemma \ref{lm 4.1} again we  that
\[\Big|\nabla\Phi_c\ast  [Q |{\bf v}+k{\bf w}_\sigma|^{p-1}({\bf v}+k{\bf w}_\sigma)]\Big|  \leq 
 \left\{\begin{array}{lll}
c_{34}k^p|x|^{\theta_p-1} \ \  &{\rm for }\ \ 0<|x|\leq 1,\\[2mm]
 \phantom{   }
c_{34}k^p |x|^{ -\sigma_p }\ \ &{\rm for }\ \ |x|> 1.
\end{array}\right.\]

So we have that $\mathcal{T}_c {\bf v}\in W^{1.\infty}_{loc}(\R^N\setminus\{0\},\C)$,
by the embeddings $W^{1,\infty}(B_{R}\setminus B_{\frac1R},\C)\hookrightarrow   L^1( B_{R}\setminus B_{\frac1R},\C)$   are  compact and
together with upper bound in $\mathcal{D}_{p,k}$, we obtain
  that $\mathcal{T}_c$ is a compact operator in $\mathcal{D}_{p,k}$.

Observing that $\mathcal{D}_{p,k}$ is a closed and convex set in $L^1(\R^N,\C)$, we may apply Schauder
fixed point theorem to derive that there exists ${\bf v}_k\in \mathcal{D}_{p,k}$ such that
$$\mathcal{T}_c{\bf v}_k={\bf v}_k. $$
Since $   |{\bf v}_k|\le  k W_p$, so $v_k$ is locally bounded in $\R^N\setminus\{0\}$,
then ${\bf u}_k:={\bf v}_k+k{\bf w}_\sigma$ satisfies 
$$\limsup_{|x|\to+\infty} |{\bf u}_k(x)||x|^{\sigma}<+\infty.$$

Moreover, if $$\alpha+\sigma p-\frac{1+N}2>    \sigma,$$
then (\ref{eq 4.1-1-1-c}) implies that 
$$|{\bf v}_k(x)|=|\Phi_c\ast [Q |{\bf v}_k+k{\bf w}_\sigma|^{p-1}({\bf v}_k+k{\bf u}_\sigma)]|\leq c_{36}k^p|x|^{-\sigma_p}\quad {\rm for}\ \, |x|>1$$
and  then
$$\limsup_{|x|\to+\infty} |{\bf u}_k(x)-k{\bf w}_\sigma(x)||x|^{\sigma_p}\in[0,+\infty).$$
By the standard interior regularity results for system (\ref{eq 1.1 comw}), ${\bf u}_k$, is a positive classical solution of (\ref{eq 1.1}).
From Theorem \ref{teo 1} it implies that ${\bf u}_k$ is a distributional solution of (\ref{eq 1.1}). \hfill$\Box$

  \bigskip  \bigskip

   \noindent{\bf Acknowledgements:} This work is is supported by NNSF of China, No: 12071189 and 12001252,
by the Jiangxi Provincial Natural Science Foundation, No: 20202BAB201005, 20202ACBL201001. F.Z. is supported by Science and
Technology Commission of Shanghai Municipality (STCSM), grant No. 18dz2271000 and also supported by
NSFC (No. 11431005).

\bibliographystyle{amsplain}

\begin{thebibliography}{99}

\bibitem{Av} P. Aviles, Local behaviour of the solutions of some elliptic equations,
{\it Comm. Math.
Phys. 108}, 177-192 (1987).


\bibitem{A} S. Agmon,  A representation theorem for solutions of the Helmholtz equation and resolvent estimates for the Laplacian, {\it  Analysis, et cetera,}  39-76 (1990).

\bibitem {BB}
Ph. B\'{e}nilan and H. Brezis, Nonlinear problems related to the
Thomas-Fermi equation, {\it J. Evolution Eq. 3}, 673-770 (2003).

\bibitem {Br} H. Brezis, Some variational problems of the Thomas-Fermi type.
Variational inequalities and complementarity problems, {\it Proc.
Internat. School, Erice, Wiley, Chichester}, 53-73 (1980).

\bibitem {BL} H. Brezis and P.L. Lions; A note on isolated singularities for linear elliptic equations, {\it in Mathematical Analysis and Applications, Acad. Press,  } 263-266 (1981).

\bibitem {BV}
H. Brezis and L. V\'{e}ron, Removable singularities of some nonlinear elliptic equations,
 {\it Arch. Rational Mech. Anal. 75}, 1-6 (1980).
 
 \bibitem{CGS} L. Caffarelli, B. Gidas and J. Spruck, Asymptotic symmetry and local behaviour of
semilinear elliptic equations with critical Sobolev growth, {\it Comm. Pure Appl. Math.
42}, 271-297 (1989).

\bibitem{CQ} H. Chen, A. Quaas, Classification of isolated singularities of nonnegative solutions to fractional semi-linear elliptic equations and the existence results,
{\it  J. London Math. Soc. 97(2)}: 196-221 (2018).

\bibitem{CEW} H. Chen, G. Ev\'{e}quoz and T. Weth,   Complex solutions and stationary scattering for the nonlinear Helmholtz equation, {\it SIAM J. Math. Anal. 53(2),} 2349-2372 (2021).



\bibitem{EP0} A. Enciso, D. Peralta-Salas, Bounded solutions to the Allen-Cahn equation with
level sets of any compact topology, {\it Analysis and PDE 9(6)}, 1433-1446 (2016).



\bibitem{E} G. Ev\'{e}quoz, A dual approach in Orlicz spaces for the nonlinear Helmholtz equation,    {\it  Z. Angew. math. Phys. 66}, 2995-3015 (2015).

\bibitem{EW0} G. Ev\'{e}quoz and T. Weth, Real solutions to the nonlinear Helmholtz equation with local nonlinearity, {\it  Arch. Rational Meth. Anal. 211}, 359-388 (2014).

\bibitem{EW1} G. Ev\'{e}quoz and T. Weth,  Dual variational methods and nonvanishing for the nonlinear Helmholtz equation, {\it   Adv. math. 280}, 690-728   (2015).


\bibitem{EY} G. Ev\'{e}quoz and T. Yesil,  Dual ground state solutions for the critical nonlinear Helmholtz equation, {\it    arXiv: 1707.00959} (2017).




\bibitem{GS} B. Gidas and J. Spruck, Global and local behaviour of positive solutions of nonlinear
elliptic equations, {\it Comm. Pure Appl. Math. 34}, 525-598 (1981).


\bibitem{G} S. Guiti\'{e}rrez, Non trivial $L^q$ solutions to the Ginzburg-Landau equation,
{\it Math. Ann. 328}, 1-25 (2004).

\bibitem{J} E. Jalade, Inverse problem for a nonlinear Helmholtz equation, 
{\it Ann. l'IHP Anal. non lin\'eaire. 21(4)},  517-531 (2004).

\bibitem {KR}  T. Klimsiak,  A. Rozkosz,   On semilinear elliptic equations with diffuse measures. {\it NoDEA Nonlinear Differential Equations Appl. 25(4)}  No. 35, 23 pp.   (2018).


\bibitem {L}  P. Lions,   Isolated singularities in semilinear problems, {\it J. Diff. Eq. 38(3)}, 441-450 (1980).



 \bibitem{lebedev} N. Lebedev,  Special functions and their applications, {\it Dover Publications Inc., New York, 1972.}




\bibitem{M1} R. Mandel, The limiting absorption principle for periodic differential operators and applications to nonlinear Helmholtz equations,
{\it Comm.   Math. Phys.  368(2)}, 799-842(2019). 

\bibitem{MMP} R. Mandel, E. Montefusco and B. Pellacci, 
Oscillating solutions for nonlinear Helmholtz equations,
{\it Z. Angew. Math. Phys. }  68-121(2017)



\bibitem {PW} A. C. Ponce and N. Wilmet,  Schr\"odinger operators involving singular potentials and measure data, J. Diff. Eq. 263, 3581-3610 (2017).


\bibitem{S} E. M. Stein, Singular Integrals and Differentiability Properties of Functions, {\it Princeton University Press,} (1970).

\bibitem{W} Y. Wang,  Isolated singularities of solutions of defocusing Hartree equation. {\it Nonlinear Anal. 156}, 70-81 (2017).


\bibitem {V1} L. V\'{e}ron, Weak and strong singularities of nonlinear elliptic equations, {\it Proc. Symp. Pure Math.
45,} 477-495 (1986).

\bibitem {V0}  L. V\'{e}ron, Singular solutions of some nonlinear elliptic equations,
{\it Nonlinear Anal. T. M. $\&$ A. 5,} 225-242 (1981).

\bibitem {V}  L. V\'{e}ron, Elliptic equations involving Measures,
 Stationary Partial Differential equations,
{\it Vol. I, 593-712, Handb. Differ. Equ., North-Holland, Amsterdam}
\end{thebibliography}

 \end{document}